\documentclass[12pt,leqno]{article}
\usepackage{amstext}
\usepackage{theorem}
\usepackage{amssymb}
\usepackage{latexsym}
\newfont{\ninemsbm}{msbm10 scaled 0900}
\newfont{\tenmsbm}{msbm10 scaled 1100}
\newfont{\nineeufb}{eufb10 scaled 0900}
\newfont{\teneufb}{eufb10 scaled 1100}
\newfont{\teneusm}{eusm10 scaled 1100}
\newfont{\nineeusm}{eusm10 scaled 0900}
\newfont{\paren}{cmex10 scaled 1100}
\newcommand{\be}{\begin{enumerate}}
\newcommand{\ee}{\end{enumerate}}

\newcommand{\br}{\begin{array}}
\newcommand{\er}{\end{array}}

\newcommand{\gchi}{\raisebox{.4ex}{$\chi$}}

\newcommand{\bdquote}{``}
\newcommand{\edquote}{''}

\newtheorem{thm}{Theorem.}[section]
\newtheorem{prop}[thm]{Proposition.}

\newtheorem{defn}[thm]{Definition.}

\def\pmb#1{\setbox0=\hbox{#1}%
\kern-.025em\copy0\kern-\wd0
\kern.05em\copy0\kern-\wd0
\kern-.025em\raise.0433em\box0}

\begin{document}

\title{Interpolation on Jets}
\author{J. Alexander \and A. Hirschowitz}
\date{}
\maketitle
\section{Introduction}

We work over a field $k$ of characteristic zero.

A {\em jet} in $n$-dimensional projective space
$\mbox{{\tenmsbm P}}_{k}^{\, n}$ over $k$
will be any divisor on a one dimensional linear 
subspace (i.e. a line in $\mbox{{\tenmsbm P}}_{k}^{\,
n}$) with support a point. The {\em length} of a jet 
will be its degree as a divisor and we will say
$r$-{\em jet} for a jet of length $r$. When $1< r$,
an $r$-jet is contained in a unique  line, which we
call its  {\em axis}. As for $r = 1$,  we will
distinguish $1$-jets, which are furnished with an
axis (namely a line containing them) from (free)
points. Finally,  by convention, a $0$-jet is the
empty subscheme.

We say that a closed sub-scheme $Y\hookrightarrow
\mbox{{\tenmsbm P}}_{k}^{n}$ has {\em maximal rank in
degree} $d\geq 0$ if the canonical map
$H^{0}(\mbox{{\tenmsbm P}}_{k}^{\, n},\mbox{{\teneusm
O}}_{\mbox{{\tenmsbm P}}^{\, n}}(d))
\longrightarrow H^{0}(Y,\mbox{{\teneusm O}}_{Y}(d))$
has maximal rank as a linear map. 

Our theorem 1.2 slightly refines a  result  proved by
A. Eastwood in [E1,2] characterizing those generic
unions of jets having maximal rank. The proof in [1,
2] used complex specialisation arguments and his
proof runs to some sixty pages. Our proof is
independent of this  earlier one and much simpler.
It works by an induction arguement based on the
simple proposition 2.1, which allows us to reduce to
the extremal case where the sequence of lengths is
maximal (among the allowed ones) with respect to the
lexicographical order. This case can then be treated
by elementary techniques when $n=2$ and for $n\geq
3$, using,  as was done in  [1, 2], using an old
theorem of Hartshorne and the second author [3].
This theorem says that generic unions of lines have
maximal rank in any degree.
\begin{defn}We will say that a set $S$ of lines is
{\em in general position} if
 for any subset of $S$, the corresponding union
has maximal rank in any degree.\end{defn}
 Observe that in $\mbox{{\tenmsbm P}}^{\,
2}$, any finite set of lines is in general position,
while in  $\mbox{{\tenmsbm P}}^{\, n}$, it follows
from the above mentioned result that any generic
finite set of lines is in general position.

\begin{thm}
Let $L = (L_{1},\ldots , L_{m})$ be a sequence of
$m$ lines in general position in $\mbox{{\tenmsbm
P}}_{k}^{\, n}$, $r = (r_1, ..., r_m)$ be  a sequence
of positive integers in  non-increasing order and
$d$ an integer such that the sum $r_{1}+\cdots +
r_{m}$ is at most 
$\scriptsize{\left(\begin{array}{c} n+d \\
d\end{array}\right)}$. Then the union $J_1 \cup ...
\cup J_m$ where $J_i$ is the generic $r_i$-jet on
$L_i$  has maximal rank in degree $d$ if and only if
the following (necessary) numerical condition
$C(n,d)$ holds: $r_{1}\leq d+1$ and, if $n=2$, then
for any $1\leq s \leq d+1$,
$$\mbox{$r_{1}+\cdots + r_{s}\leq d\, .\,
s+1-$}\scriptsize{\left(\begin{array}{c} s-1
\\ 2
\end{array}\right)}.$$
\end{thm}

In other words, if one fixes the degree
$d$, then for given lines $L_1,\ldots ,L_m\in
\mbox{{\tenmsbm P}}^n$ in general position, the
$\Sigma_i r_i$ conditions imposed on hypersurfaces
of degree $d$ by requiring a contact of order
$r_i$ at the generic point $x_i$ of $L_i$ for
$i=1,\ldots ,m$, are linearly independent precisely
under the expected numerical conditions given in the
theorem.

Before passing on to the proof,
we introduce the following notation.
\begin{defn} For a fixed $n$, we
say that a sequence of non-negative integers $(\gchi
,r_1,\ldots ,r_m)$ is $d$-{\em admissible} if
\newcounter{mist0}
\begin{list}{\textbf{ (\roman{mist0}) }}
{\setlength{\leftmargin}{1cm}
\setlength{\rightmargin}{1cm}
\usecounter{mist0}  } 
\item 
$\mbox{$\gchi + r_1\cdots
+
r_m$}=\scriptsize{\left(\begin{array}{c}n+d
\\d\end{array}\right)}$
\item $r_1\geq \cdots \geq r_m \geq 1$
\item if $n=2$ then
$$r_1+\cdots r_s\leq ds+1-(s-1)(s-2)/2$$
for $s=1,\ldots ,d+1$.
\end{list}
and we denote by $S_d$ the set of all $d$-admissible
sequences, which we give the total lexicographical
order (i.e.  $(\alpha_i)_i> (\beta_i)_i$ if the
first non-zero term in the sequence
$(\alpha_i-\beta_i)_i$ is strictly positive).
\end{defn}

\begin{defn} Given a $d$-admissible sequence
$\mbox{{\boldmath $r$}}=(\gchi , r_1,\ldots ,r_m)$
and a sequence of lines in general position
$L = (L_{1},\ldots , L_{m})$ it is clear that the
union $J_1 \cup \cdots \cup J_m$ of the generic
$r_i$-jets $J_i$ in $L_i$, has maximal rank in
degree $d$ if and only if the union $J :=
J_{L,r}=J_1 \cup \cdots. \cup J_m \cup R$  has
maximal rank, where $R$ is the union  of $\gchi$
generic (free) points. We call such a union
$J$ {\em a
$d$-admissible union of jets} and we say that 
$$\mbox{{\boldmath $r$}}=(\gchi , r_1,\ldots ,r_m)$$
is the weight of $J$. 
\end{defn}

In this language, theorem 1.2
simply says that every $d$-admissible union of jets
has maximal rank in degree $d$.

Without further comment, we will freely use the fact
that \bdquote maximal rank in degree
$d$ \edquote is stable by generisation.
We express our gratitude to the referee for having
suggested changes which have improved the
presentation of the proof.
\section{The proof}

\begin{prop} Fix a closed subscheme $Y\subset
\mbox{{\tenmsbm P}}_{k}^{\, n}$ and
let $D$ be a line not contained in $Y$.  Suppose
that for some $d, v >0$ the union of $Y$ with the
generic jet of length $v+1$ (resp $v-1$) in $D$
 has maximal rank in degree $d$,
then the union of $Y$ with the generic $v$-jet in
$D$ has maximal rank in degree
$d$.
\end{prop}

\noindent{\bf Proof.}
We denote by $Y_r$ the union of $Y$ with the generic
$r$-jet; denoted $J_{r}$; in $D$ and by
$m_r$ the canonical map
$$H^0(\mbox{{\tenmsbm P}}^n,\mbox{{\teneusm
O}}_{\mbox{{\ninemsbm P}}^n}(d))
\longrightarrow H^0(\mbox{{\tenmsbm
P}}^n,\mbox{{\teneusm O}}_{Y}(d))\times
H^0(\mbox{{\tenmsbm P}}^n,\mbox{{\teneusm
O}}_{J_{r}}(d))\, =\, 
 H^0(\mbox{{\tenmsbm P}}^n,\mbox{{\teneusm
O}}_{Y_r}(d))$$
If $m_{v+1}$ is surjective, then so is $m_v$, and if
$m_{v-1}$ is injective, then so is $m_v$. In the
remaining case, $m_{v+1}$ is injective and not
surjective,  and  $m_{v-1}$ is surjective and not
injective and we have to prove that $m_v$ is
bijective. 

What is clear is that $h^0(\mbox{{\tenmsbm P}}^n,
I_Y(d))=v$, where $I_Y$ is the ideal sheaf of $Y$ as
a subscheme of $\mbox{{\tenmsbm P}}^n$, and that the
canonical map
$$H^0(\mbox{{\tenmsbm P}}^n, I_Y(d))\hookrightarrow
H^0(\mbox{{\tenmsbm P}}^n,
\mbox{{\teneusm O}}_{J_{v+1}}(d))$$
is injective. It follows that the map
$$H^0(\mbox{{\tenmsbm P}}^n, I_Y(d))\hookrightarrow
H^0(D,
\mbox{{\teneusm O}}_{D}(d))$$
is injective. If $W$ is the image subspace of this
latter map we need only show that the generic
$v$-jet in $D$ imposes linearly independent
conditions on $W$, but this is true for any $v$
dimensional subspace of $H^0(D,\mbox{{\teneusm
O}}_{D}(d))$. This follows from the well known fact
that a familly of polynomials
$f_1,\ldots ,f_v\in k[t]$ is linearly independent if
and only if the $v\times v$ Wronskian 
$$\begin{array}{lll}W(f_1,\ldots
,f_v)&=&\left[\begin{array}{l}\partial^i f_j
\vspace{.5ex}\\ \hline\vspace{-2ex} \\ \partial t^i
\end{array}\right]\end{array}$$  has rank $v$.

\subsection{Proof of theorem 1.2 for
$\mbox{{\tenmsbm P}}^{2}$}

As indicated after definition 1.4, we will prove
that every $d$-admissible  union  of jets  has
maximal rank in degree $d$.

We will write a $d$-admissible weight
$\mbox{{\boldmath
$r$}}=(\gchi,\sigma_1,\ldots ,\sigma_m)$ in the form
$$\mbox{{\boldmath $r$}}=(\gchi,m_1,\ldots
,m_p,r_1,\ldots ,r_q)$$ where $m_1,\ldots ,m_p$ is
the extremal sequence $d+1,\ldots ,d+2-p$ given by
the condition 1.3 (iii) and $r_1< d+1-p$.

Fix a $d$-admissible union of jets $J$ of weight
$$\mbox{\boldmath $r$}=(\gchi,m_1,\ldots ,
m_p,r_1,\ldots , r_q)$$ and suppose that every
$d$-admissible union of jets $J^{\prime}$ of weight  
$\mbox{{\boldmath $r$}}^{\prime}>\mbox{{\boldmath
$r$}}$ has maximal rank in degree $d$. We will show
that this implies that $J$ has maximal rank in
degree $d$ completing the proof.

If $q\leq 1$ (i.e. $r_2=0$) then by specialising the
$\gchi$ free points, we can specialise $J$ to the
$d$-admissible union of jets
$J_d$ of weight $(0,d+1,\ldots ,1)$. A simple
induction then shows that any homogeneous form of
degree
$d$ on $\mbox{{\tenmsbm P}}^2$ which vanishes on
$J_d$ vanishes on all $d+1$ axes and is thus
identically zero. This shows that $J$ has maximal
rank in degree
$d$ if $r_2=0$.

Now suppose that $r_2>0$. We will apply 2.1 with
$$Y=J_1\cup\cdots \cup J_p\cup
J^{\prime\prime}_2\cup \cdots \cup
J_q^{\prime\prime}\cup R$$ where $J_i$ is the
$m_i$-jet, $J_i^{\prime\prime}$ is the $r_i$-jet and
$R$ is the union of the $\gchi$ free points. We must
show that the unions
$J_{(+)}$ and $J_{(-)}$ of $Y$ with the jet of
length $r_1+1$ (resp. $r_1-1$) extending (resp.
contracting) $J_1^{\prime\prime}$ have maximal rank
in degree
$d$. 

On the one hand, $J_{(-)}$ is contained in the union
$\widetilde{J}_{(-)}$ obtained by adding a generic
free point to $J_{(-)}$. Now $\widetilde{J}_{(-)}$
is clearly a $d$-admissible union of jets of weight 
$\widetilde{\mbox{{\boldmath$r$}}}_{(-)}>
\mbox{{\boldmath$r$}}$
so that
$J_{(-)}$ has maximal rank in degree $d$ by the
induction hypothesis.

On the other hand, $J_{(+)}$ contains the union of
jets $\widetilde{J}_{(+)}$ obtained by contracting
$J_2^{\prime\prime}$ to length $r_2-1$ and it will
suffice to show that the weight
$$\widetilde{\mbox{{\boldmath$r$}}}_{(+)}
=(\gchi,m_1,\ldots ,m_p,r_1+1,r_2,\ldots
,r_{t-1},r_t-1,r_{t+1},\ldots ,r_q)$$ of
$\widetilde{J}_{(+)}$ is $d$-admissible, where
$r_2=\cdots = r_t>r_{t+1}$. This is
equivalent to showing that the series of inequalities
$$r_1 + (s-1)r_2 \leq (d-p)s+1-(s-1)(s-2)/2$$
for $1\leq s\leq \mbox{min}\, (t,d-p+1)$ implies
strict inequality for 
$1\leq s\leq \mbox{min}\, (t-1,d-p+1)$. Since the
associated quadratic form
$$Q(s)=s^2 -(3+2(d-p-r_2)) s + 2(r_1-r_2)$$
vanishes on the interval $[ 0,1[$, there is only one
root $\geq 1$. This completes the proof if $t\leq
d-p+1$ or if $Q$ has no integer roots $\leq d-p+1$.
We will now finish by showing that if $t> d-p+1$
then $Q$ does not have integer roots. In fact if $Q$
has integer roots, they must be
$0$ and $s_1=3+2(d-p-r_2)\leq d-p+1$. In this case
$r_2\geq (d-p+2)/2$. However we have 
$$tr_2\leq r_1+(t-1)r_2\leq r_1+\ldots + r_t\leq
(d-p+2)(d-p+1)/2$$ showing that $t\leq
d-p+1$.$\;\;\;\;\;\;\;\Box$

\subsection{Proof of theorem 1 for $\mbox{{\tenmsbm
P}}^{n}$, $n\geq 3$.} The proof is similar to the
previous one. We write the weight 
$\mbox{{\boldmath $r$}}$ of a $d$-admissible
union of jets $J$ in the form
$$\mbox{{\boldmath $r$}}=(\gchi, m_1,\ldots ,
m_p,r_1,\ldots , r_q)$$ where $m_i=d+1$ and $r_1\leq
d$.

Now let $J$ be  a $d$-admissible
union of jets of weight $\mbox{{\boldmath $r$}}$ and
suppose that every
$d$-admissible union of jets $J^{\prime}$ of weight
$\mbox{{\boldmath
$r$}}^{\prime}> \mbox{{\boldmath $r$}}$ has maximal
rank in degree $d$. We will show that this implies
that $J$ has maximal rank in degree $d$.

If $r_2\neq 0$ one easily applies 2.1 as before. So
suppose that $r_2=0$. By specialising the $\gchi$
free points firstly to the $r_1$-jet then to further
lines in general position we can specialise $J$ to a
$d$-admissible union of jets of weight 
$\mbox{{\boldmath $r$}}=(0,m_1,\ldots ,m_p,\delta)$
where
$0\leq \delta\leq d$. This is equivalent to the
corresponding problem obtained by replacing  the
$p$, $(d+1)$-jets by their corresponding axes. That
is to say, we must show that if $D_1,\ldots,D_{p+1}$
are $p+1$ lines in general position, then the union
of $D_1,\ldots ,D_p$ with the generic $\delta$-jet
in $D_{p+1}$ has maximal rank in degree $d$. The
general position hypothesis implies that 
$h^0(\mbox{{\tenmsbm
P}}^n,I_{D_1\cup\, \cdots \, \cup D_p}(d))\; =\;
\delta$ and that the canonical map
$$H^0(\mbox{{\tenmsbm
P}}^n,I_{D_1\cup\, \cdots\, \cup
D_p}(d))\longrightarrow H^0(D_{p+1}, \mbox{{\teneusm
O}}_{D_{p+1}}(d))$$ is injectif. This completes the
proof as in the concluding arguement of the proof of
proposition 2.1. $\;\;\;\;\;\;\;\;\Box$\vspace{2ex}\\

Note: our method would work equally well for any
union of curves: from a maximal rank statement
concerning a union of curves, we derive a maximal
rank statement for unions of generic jets on these
curves, under natural necessary numerical conditions.

\noindent{\Large\bf Bibliography.}\\
\begin{enumerate}
\item  Eastwood  A.: {\em Collision de biais et
interpolation} Manuscr. Math. 67 (1990), 227-249.
\item  Eastwood A. : {\em Interpolation \`a N
variables } J. of Algebra, 139 (1991) 273-310.
\item  Hartshorne R., Hirschowitz  A. : {\em Droites
en position g\'en\'erale}, in Proceedings La Rabida
1981, LNM 961, 169-189.
\item Hirschowitz A. : {\em Probl\`emes de
Brill-Noether en rang sup\'erieur} Pr\'eprint Nice
1986.
\end{enumerate}
$\;$\vspace{1cm}\\

\begin{minipage}[t]{.4\linewidth}{James
Alexander\\
Universit\'e d'Angers\\
jea@tonton.univ-angers.fr}\end{minipage}\hfill
\begin{minipage}[t]{.4\linewidth}{Andr\'e
Hirschowitz\\
Universit\'e de Nice\\
ah@math.unice.fr}\end{minipage}

\end{document}